\documentclass{birkmult}
 \newtheorem{thm}{Theorem}[section]

 \theoremstyle{definition}
 
 \theoremstyle{remark}

 \numberwithin{equation}{section}
 \def\RR{{\mathbb R} }

\begin{document}
%
%
%
%
%
%
%
%
%
\title[Bifurcation and Asymptotics for Singular Elliptic Problems]
 {Bifurcation and Asymptotics for Elliptic\\ Problems with
Singular Nonlinearity}
\author[Vicen\c tiu R\u adulescu]{Vicen\c tiu R\u adulescu}

\address{%
Department of Mathematics\\University of Craiova\\RO 200585
Craiova\\Romania\\{\tt http://inf.ucv.ro/\~{}radulescu}}

\email{radulescu@inf.ucv.ro}

\thanks{Partially supported by Grant 12/2004 with the Romanian Academy.}
\subjclass{Primary 35J60; Secondary 35B32, 35B40}

\keywords{Singular nonlinearity, bifurcation, asymptotic analysis,
maximum principle.}

\date{October 1, 2004}
\dedicatory{A mon Ma\^{\i}tre, avec reconnaissance}

\begin{abstract}
We report on some recent existence and uniqueness results for
elliptic equations subject to Dirichlet boundary condition and
involving a singular nonlinearity. We take into account the
following types of problems: (i) singular problems with sublinear
nonlinearity and two parameters; (ii) combined effects of
asymptotically linear and singular nonlinearities in bifurcation
problems; (iii) bifurcation for a class of singular elliptic
problems with subquadratic convection term. In some concrete
situations we also establish the asymptotic behaviour of the
solution around the bifurcation point. Our analysis relies on the
maximum principle for elliptic equations combined with adequate
estimates.
\end{abstract}

\maketitle
\section{Motivation and Previous Results}
I will report on some  results contained in our recent papers
\cite{cgr,gr1,gr2,gr3, gr4,gr5} that are closely related to the
study of some problems on blow-up boundary solutions. More
precisely, consider the elementary example
$$\left\{\begin{tabular}{ll}
$\Delta u=u^p$ \quad & $\mbox{\rm in}\ \Omega,\ $\\
$u>0$ \quad & $\mbox{\rm in}\ \Omega,$\\
$u=+\infty$ \quad & $\mbox{\rm on}\ \partial\Omega\,,$
\end{tabular} \right.$$
where $\Omega\subset\RR^N$ is a smooth bounded domain and $p>1$.
Then the function $v=u^{-1}$ satisfies
\begin{equation}\label{eq1}\left\{\begin{tabular}{ll}
$-\Delta v=v^{2-p}-\displaystyle  \frac 2v\,|\nabla v|^2$ \quad & $\mbox{\rm in}\ \Omega,$\\
$v>0$ \quad & $\mbox{\rm in}\ \Omega,$\\
$v=0$ \quad & $\mbox{\rm on}\ \partial\Omega.$\\
\end{tabular} \right.\end{equation}
The above equation contains both singular nonlinearities (like
$v^{-1}$ or $v^{2-p}$, if $p>2$) and a convection term (denoted by
$|\nabla v|^2$). These nonlinearities make more difficult to
handle problems like \eqref{eq1}. Our purpose in this paper is to
give an overview on some old and new results in this direction. We
recall the pioneering paper \cite{crt} that contains one of the
first existence results for singular elliptic problems. In fact,
it is proved in \cite{crt} that the boundary value problem
$$
 \left\{\begin{tabular}{ll}
$-\Delta u-u^{-\alpha}= -u$ \quad & ${\rm in}\
\Omega,$\\
$u>0$ \quad & ${\rm in}\ \Omega,$\\
$u=0$ \quad & ${\rm on}\ \partial\Omega$
\end{tabular} \right.
$$
has a solution, for any $\alpha>0$. Let us now consider the
problem
\begin{equation} \label{unuunu}
 \left\{\begin{tabular}{ll}
$-\Delta u-u^{-\alpha}=\lambda u^p$ \quad & ${\rm in}\
\Omega,$\\
$u>0$ \quad & ${\rm in}\ \Omega,$\\
$u=0$ \quad & ${\rm on}\ \partial\Omega,$\\
\end{tabular} \right.
\end{equation}
where $\lambda\geq 0$ and $\alpha,p\in(0,1).$ In \cite{cp} it is
proved that problem \eqref{unuunu} has at least one solution for
all $\lambda\geq 0$ and $0<p<1$. Moreover, if $p\geq 1,$ then
there exists $\lambda^*$ such that problem \eqref{unuunu} has a
solution for $\lambda\in[0,\lambda^*)$ and no solution for
$\lambda>\lambda^*$. In \cite{cp} it is also proved a related
non--existence result. More exactly, the  problem
$$
 \left\{\begin{tabular}{ll}
$-\Delta u+u^{-\alpha}= u$ \quad & ${\rm in}\
\Omega,$\\
$u>0$ \quad & ${\rm in}\ \Omega,$\\
$u=0$ \quad & ${\rm on}\ \partial\Omega$\\
\end{tabular} \right.
$$
has no solution, provided that $0<\alpha <1$ and $\lambda_1\geq 1$
(that is, if $\Omega$ is ``small''), where $\lambda_1$ denotes the
first eigenvalue of $(-\Delta)$ in $H^1_0(\Omega)$.

Problems related to multiplicity and uniqueness become difficult
even in simple cases. In \cite{shi} it is studied the existence of
radial symmetric solutions to the problem
$$
 \left\{\begin{tabular}{ll}
$\displaystyle \Delta u+\lambda(u^p-u^{-\alpha})=0$ \quad & ${\rm
in}\ B_1,$\\
$u>0$ \quad & ${\rm in}\ B_1,$\\
$u=0$ \quad & ${\rm on}\ \partial B_1,$\\
\end{tabular} \right.
$$
where $\alpha>0$, $0<p<1,$ $\lambda>0$, and $B_1$ is the unit ball
in $\RR^N$. Using a bifurcation theorem of Crandall and
Rabinowitz, it has been shown in \cite{shi}  that there exists
$\lambda_1>\lambda_0>0$ such that the above problem has no
solutions for $\lambda<\lambda_0,$ exactly one solution for
$\lambda=\lambda_0$ or $\lambda>\lambda_1,$  and two solutions for
$\lambda_0 < \lambda\leq\lambda_1$.

Our purpose in this survey paper is to present various existence,
 and non--existence results  for several
classes of singular
 elliptic problems. We also take into account bifurcation nonlinear problems
 and establish the precise rate decay of the solution
 in some concrete situations. We intend to reflect the
``competition" between different quantities, such as: sublinear or
superlinear nonlinearities, singular nonlinear terms  (like
$u^{-\alpha}$, for $\alpha>0$), convection nonlinearities  (like
$|\nabla u|^{q}$, with $0<q\leq 2$), as well as sign--changing
potentials.

\section{A Singular Problem with Sublinear Nonlinearity}
Consider the following boundary value problem with two parameters:
\begin{equation}\label{Plamu}
\left\{\begin{tabular}{ll} $-\Delta
u+K(x)g(u)=\lambda f(x,u)+\mu h(x)$ \quad & ${\rm
in}\ \Omega,$\\
$u>0$ \quad & ${\rm in}\ \Omega,$\\
$u=0$ \quad & ${\rm on}\ \partial\Omega,$\\
\end{tabular} \right. \end{equation}
where $\Omega$ is a smooth bounded domain in $\RR^N$ ($N\geq 2$),
$K,h\in C^{0,\gamma}(\overline\Omega ),$ with $h>0$ on $\Omega$,
and $\lambda,\,\mu$ are positive real numbers. We suppose that
$\,f:\overline{\Omega}\times[0,\infty)\rightarrow[0,\infty)\,$ is
a H\"{o}lder continuous function which is positive on
$\,\overline{\Omega}\times(0,\infty)$. We also assume that $\,f\,$
is non--decreasing with respect to the second variable and is
sublinear, that is,

\smallskip
\noindent $\displaystyle (f1)\qquad$ the mapping $\displaystyle
(0,\infty)\ni s\longmapsto\frac{f(x,s)}{s}\quad\mbox{is
non--increasing for all}\;\, x\in\overline{\Omega};$

\smallskip \noindent $\displaystyle (f2)\qquad \lim_{s\downarrow
0}\frac{f(x,s)}{s}=+\infty\quad\mbox{and}\;\;
\lim_{s\rightarrow\infty}\frac{f(x,s)}{s}=0,\;\;\mbox{uniformly
for}\;\,x\in\overline{\Omega}.$

\smallskip We assume that $g\in C^{0,\gamma}(0,\infty)$ is a
non--negative and non--increasing function. A fundamental role in
our analysis will be played by the numbers
$$K^*:=\max_{x\in{\overline{\Omega}}}K(x),\qquad K_*=\min_{x\in{\overline{\Omega}}}K(x).$$

Our first theorem is a non--existence result and it concerns
nonlinearities with strong blow-up rate at the origin (like
$u^{-\alpha}$, with $\alpha\geq 1$).

\begin{thm}\label{th1}
Assume that $K_*>0$ and $f$ satisfies $(f1)-(f2)$. If $\int^1_0
g(s)ds=+\infty,$ then problem \eqref{Plamu} has no classical
solution, for any $\lambda$, $\mu>0.$
\end{thm}

Next, we assume that the growth of the nonlinearity is described
by the following conditions:

\smallskip
 \noindent $\displaystyle (g1)\qquad
\lim_{s\downarrow 0}g(s)=+\infty;$

\smallskip \noindent $\displaystyle (g2)\qquad \mbox{there
exist }C,\ \delta_0>0$ and $\alpha\in(0,1)$ such that $g(s)\leq
Cs^{-\alpha},$ for all $s\in(0,\delta_0).$

\smallskip The above
conditions $(g1)$ and $(g2)$ are fulfilled by singular
nonlinearities like $g(u)=u^{-\alpha}$, with $\alpha\in (0,1)$.
Obviously, hypothesis $(g2)$ implies
 the following Keller-Osserman type condition around
the origin:

\smallskip \noindent $\displaystyle (g3)\qquad \int^1_{0}\left(
\int^t_0g(s)ds\right)^{-1/2}dt<\infty.$

 As proved by B\'enilan, Brezis and Crandall \cite{bbc},
condition $(g3)$ is equivalent to the {\it property of compact
support}, that is, for every $h\in L^1(\RR^N)$ with compact
support, there exists a unique $u\in W^{1,1}(\RR^N)$ with compact
support such that $\Delta u\in L^1(\RR^N)$ and $-\Delta u+g(u)=h$,
a.e. in $\RR^N.$ That it is why it is natural to try to find
solutions in the class $${\mathcal E}=\{\,u\in C^2(\Omega)\cap
C({\overline{\Omega}});\,\,\Delta u\in L^1(\Omega)\}.$$

In the case where the potential $K(x)$ has a constant sign, the
following results hold.

\begin{thm}\label{th2}
Assume that $ K_*>0,$ $ f $ satisfies $ (f1)-(f2) $, and $ g $
satisfies $ (g1)-(g2).$ Then there exists $ \lambda_*,\mu_*>0 $
such that:

-- problem \eqref{Plamu} has at least one solution in $ {\mathcal
E} $ either if $ \lambda>\lambda_* $ or if $ \mu>\mu_*.$

-- problem \eqref{Plamu} has no solution in $ {\mathcal E} $ if
$\lambda<\lambda_* $ and $ \mu<\mu_*.$

Moreover, if either $ \lambda>\lambda_* $ or if $ \mu>\mu_*, $
then problem \eqref{Plamu} has a maximal solution in $ {\mathcal
E} $ which is increasing with respect to $ \lambda $ and $ \mu.$
\end{thm}

\begin{figure}[h]
\begin{center}
\begin{picture}(180,180)\setlength{\unitlength}{1pt}
\put(7,-15){\makebox(40,14){(0,0)}}
\put(175,-15){\makebox(40,15){$\lambda$}}
\put(100,-15){\makebox(40,15){$\lambda_*$}} \vector(1,0){200}
\put(-190,-10){ \vector(0,20){200}}
\put(-213,175){\makebox(40,15){$\mu$}}
\put(-213,40){\makebox(40,15){$\mu_*$}}
\multiput(-186,50)(3,0){34}{.} \multiput(-85,3)(0,3){17}{.}
\put(-152,20){\makebox(40,15){No solution}}
\put(-100,80){\makebox(40,15){At least one solution}}
\end{picture}
\end{center}
\caption{The dependence on $ \lambda $ and $ \mu$ in Theorem
\ref{th2}}
\end{figure}
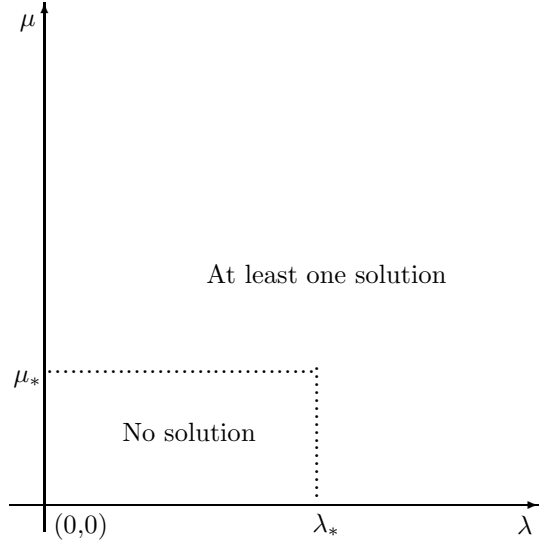

At this stage we are not able to describe the behaviour in the
following cases: (i) [$\lambda=\lambda_*$ and $0<\mu\leq\mu_*$]
and (ii)  [$0<\lambda\leq\lambda_*$ and $\mu=\mu_*$]. We
conjecture that existence or non--existence results can be
established in conjunction with a more precise description of the
decay rate of the potential coefficients and nonlinearities.

\begin{thm}\label{th3}
Assume that $ K^*\leq 0,$ $ f $ satisfies conditions $ (f1)-(f2) $
and $ g $ satisfies $ (g1)-(g2).$ Then problem \eqref{Plamu} has a
unique solution $ u_{\lambda,\mu}$ in ${\mathcal E} $, for any
$\lambda$, $\mu>0.$ Moreover, $ u_{\lambda,\mu} $ is increasing
with respect to $ \lambda $ and $ \mu.$
\end{thm}

The following result give partial answers in the case where the
potential $K(x)$ changes sign.

\begin{thm}\label{th4}
Assume that $K^*>0>K_*,$ $f$ satisfies $(f1)-(f2)$ and $g$
verifies $(g1)-(g2).$ Then there exist $\lambda_*$ and $\mu_*>0$
such that problem \eqref{Plamu} has at least one solution $
u_{\lambda,\mu}\in{\mathcal E} $, provided that either $
\lambda>\lambda_* $ or $ \mu>\mu_*.$ Moreover, for
$\lambda>\lambda_*$ or $\mu>\mu_*,$ $u_{\lambda,\mu}$ is
increasing with respect to $\lambda$ and $\mu.$
\end{thm}

The proofs of the above results rely on the sub-- and
super--solution method for elliptic equations combined with
adequate comparison principles. We refer to \cite{gr2} for
complete details and additional results.

A natural question is to see what happens if assumption $(f1)$
holds true, but if $\lim_{s\rightarrow\infty}f(x,s)/s$ is {\bf
not} zero. We give in what follows a precise description in the
case where $K\leq 0$. More exactly, we consider the problem
\begin{equation} \label{Pla}
\left\{\begin{aligned} & -\Delta u=\lambda f(u)+a(x)g(u) &&
{\rm in}\ \Omega,\\
& u>0 && {\rm in}\ \Omega,\\
& u=0 && {\rm on}\ \partial\Omega,\\
\end{aligned} \right.
\end{equation}
where $a\in C^{0,\gamma}(\overline\Omega)$, $a\geq 0$,
$a\not\equiv 0$ in $\overline\Omega$, and

\smallskip \noindent $\displaystyle (f3)\qquad
\lim_{s\rightarrow\infty}\frac{f(s)}{s}=m\in (0,\infty).$

\smallskip
Let $\lambda_1$ be the first Dirichlet eigenvalue of $(-\Delta)$
in $\Omega$ and $ \lambda^*:=\lambda_1/m$. Set
$a_*:=\min_{x\in\overline\Omega}a(x)$ and $d(x):=\mbox{dist}\,
(x,\partial\Omega)$.

\begin{thm} \label{th5} Assume that conditions $(f1)$, $(f3)$, $(g1)$, and $(g2)$
are fulfilled. Then the following hold.

\begin{enumerate}
\item[{\rm (i)}] If $\lambda\geq\lambda^*$, then problem \eqref{Pla} has no
solutions in  ${\mathcal E}$.

\item[{\rm (ii)}] If $a_*>0$ (resp.
$a_*= 0$) then problem \eqref{Pla} has a unique solution
$u_{\lambda}\in{\mathcal E}$ for all $ -\infty<\lambda<\lambda^*$
(resp. $ 0<\lambda< \lambda^*$) with the properties:

\noindent {\rm (ii1)} $u_{\lambda}$ is strictly increasing with
respect to $\lambda$;

\noindent {\rm (ii2)} there exist two positive constants $c_1$, $
c_2>0$ depending on $\lambda$ such that $c_1\,d(x)\leq
u_{\lambda}(x)\leq c_2\, d(x)$, for all  $x\in\Omega$;

\noindent {\rm (ii3)} $ \lim_{\lambda\nearrow
\lambda^*}u_{\lambda}=+\infty$, uniformly on compact subsets of
$\Omega$.
\end{enumerate}
\end{thm}

\begin{proof} The first part of the proof relies on standard
arguments based on the maximum principle (see \cite{cgr} for
details). The most interesting part of the proof concerns (ii3)
and, due to the special character of our problem, we will be able
to show that, in this case, $L^2$--boundedness implies
$H^1_0$--boundedness! We refer to \cite{mr} for a related problem
and further results.

Let $u_{\lambda}\in{\mathcal E}$ be the unique solution of
$\eqref{Pla}$ for $0<\lambda<\lambda^*$. We prove that
$\displaystyle \lim_{\lambda\nearrow
\lambda^*}u_{\lambda}=+\infty$, uniformly on compact subsets of
$\Omega$. Suppose the contrary. Since
$(u_{\lambda})_{0<\lambda<\lambda^*}$ is a sequence of nonnegative
super--harmonic functions in $\Omega$ then, by Theorem~4.1.9 in
\cite{h}, there exists a subsequence of
$(u_{\lambda})_{\lambda<\lambda^*}$ [still denoted by
$(u_{\lambda})_{\lambda<\lambda^*}$] which is convergent in
$L^1_{\rm loc}(\Omega)$.

We first  prove that $(u_{\lambda})_{\lambda<\lambda^*}$ is
bounded in $L^2(\Omega)$. We argue by contradiction. Suppose that
$(u_{\lambda})_{\lambda<\lambda^*}$ is not bounded in
$L^2(\Omega).$ Thus, passing eventually at a subsequence we have
$u_{\lambda}=M(\lambda)w_{\lambda},$ where \begin{equation}
\label{cdoi}
M(\lambda)=||u_{\lambda}||_{L^2(\Omega)}\rightarrow\infty\quad
\mbox{ as $\lambda\nearrow\lambda^*$  and $w_{\lambda}\in
L^2(\Omega),$ $\|w_{\lambda}\|_{L^2(\Omega)}=1$}.\end{equation}

Using $(f1)$, $(g2)$ and the monotonicity assumption on $g$, we
deduce the existence of $A$, $B$, $C$, $D>0$ $(A>m)$ such that
\begin{equation} \label{ctrei}
f(t)\leq At+B,\quad g(t)\leq Ct^{-\alpha}+D,\quad \mbox{ for
all}\;t>0. \end{equation} This implies
$$ \frac{1}{M(\lambda)}\left(\lambda
f(u_{\lambda})+a(x)g(u_{\lambda})\right)\rightarrow 0 \quad \mbox{
in } L^1_{\rm loc}(\Omega)\ \mbox{as }\lambda\nearrow \lambda^*$$
that is,
\begin{equation} \label{cpatru} -\Delta
w_{\lambda}\rightarrow 0 \quad\mbox{in } L^1_{\rm loc}(\Omega)\
\mbox{as }\lambda\nearrow \lambda^*.\end{equation} By Green's
first identity, we have
\begin{equation} \label{car1}
\int_\Omega \nabla w_\lambda\cdot \nabla \phi\,dx=-\int_\Omega
\phi\, \Delta w_\lambda\,dx=-\int_{ {\rm Supp}\, \phi} \phi\,
\Delta w_\lambda\,dx\quad \forall \phi\in C_0^\infty(\Omega).
\end{equation}
Using (\ref{cpatru}) we derive that
\begin{equation} \label{car2}
\begin{aligned}
\left|\int_{{\rm Supp}\,\phi} \phi\, \Delta w_\lambda\,dx\right| &
\leq
\int_{{\rm Supp}\,\phi}|\phi||\Delta w_\lambda|\,dx\\
& \leq \|\phi\|_{L^\infty}\int_{{\rm Supp}\,\phi}|\Delta
w_\lambda|\,dx\to 0\quad \mbox{as }\lambda\nearrow \lambda^*.
\end{aligned}
\end{equation}
Combining (\ref{car1}) and (\ref{car2}), we arrive at
\begin{equation} \label{car3}
\int_{\Omega} \nabla w_\lambda\cdot \nabla \phi\,dx \to 0 \ \
\mbox{as }\lambda\nearrow \lambda^*,\quad \forall \phi\in
C^\infty_0(\Omega). \end{equation} By definition, the sequence
$(w_{\lambda})_{0<\lambda<\lambda^*}$ is bounded in $L^2(\Omega)$.

We claim that $(w_{\lambda})_{\lambda<\lambda^*}$ is bounded in
$H^1_0(\Omega)$. Indeed, using \eqref{ctrei} and H\"older's
inequality, we have
$$\begin{aligned}
\int_{\Omega}|\nabla w_{\lambda}|^2&
=-\int_{\Omega}w_{\lambda}\Delta w_{\lambda}
=\frac{-1}{M(\lambda)}\int_{\Omega} w_{\lambda}
\Delta u_{\lambda}\\
& =\frac{1}{M(\lambda)}\int_{\Omega}\left[\lambda
w_{\lambda}f(u_{\lambda})+a(x)g(u_{\lambda})w_{\lambda}\right]\\
& \leq \frac{\lambda}{M(\lambda)}\int_{\Omega}
w_{\lambda}(Au_{\lambda}+B)+\frac{||a||_{\infty}}{M(\lambda)}
\int_{\Omega}w_{\lambda}(Cu^{-\alpha}_{\lambda}+D)\\
&=\lambda A\int_{\Omega}
w^2_{\lambda}+\frac{||a||_{\infty}C}{M(\lambda)^{1+\alpha}}
\int_{\Omega}w^{1-\alpha}_{\lambda}+\frac{\lambda
B+\|a\|_{\infty}D}{M(\lambda)}\int_{\Omega}w_{\lambda}\\
& \leq\lambda^*A+\frac{||a||_{\infty}C}{M(\lambda)^{1+\alpha}}
|\Omega|^{(1+\alpha)/2}+\frac{\lambda
B+\|a\|_{\infty}D}{M(\lambda)}|\Omega|^{1/2}.
\end{aligned} $$
From the above estimates, it is easy to see that
$(w_{\lambda})_{\lambda<\lambda^*}$ is bounded in $H^1_0(\Omega)$,
so the claim is proved. Then, there exists $w\in H^1_0(\Omega)$
such that (up to a subsequence) \begin{equation}\label{ccinci}
w_{\lambda}\ \rightharpoonup\ w \quad\mbox{ weakly in
}\;\;H^1_0(\Omega) \ \mbox{ as } \lambda\nearrow \lambda^*
\end{equation} and, since $H_0^1(\Omega)$ is compactly embedded
in $L^2(\Omega)$, \begin{equation}\label{csase}
w_{\lambda}\rightarrow w \quad\mbox{ strongly in }\;\;L^2(\Omega)
\ \mbox{ as } \lambda\nearrow \lambda^*.\end{equation} On the one
hand, by \eqref{cdoi} and \eqref{csase}, we derive that
$\|w\|_{L^2(\Omega)}=1$. Furthermore, using \eqref{car3} and
\eqref{ccinci}, we infer that
$$ \int_{\Omega} \nabla w\cdot \nabla \phi\,dx=0,\qquad \mbox{for all
$\phi\in C^\infty_0 (\Omega)$}. $$ Since $w\in H_0^1(\Omega)$,
using the above relation and the definition of $H_0^1(\Omega)$, we
get $w=0$. This contradiction shows that
$(u_{\lambda})_{\lambda<\lambda^*}$ is bounded in $L^2(\Omega)$.
As above for $w_\lambda$, we can derive that $u_\lambda$ is
bounded in $H_0^1(\Omega)$. So, there exists $u^*\in
H_0^1(\Omega)$ such that, up to a subsequence,
\begin{equation}\label{car4} \left\{
\begin{array}{lll} &\displaystyle u_\lambda\
\rightharpoonup\ u^*\ \ &\displaystyle\mbox{weakly in
$H_0^1(\Omega)$
as $\lambda\nearrow \lambda^*$},\\
&\displaystyle u_\lambda \rightarrow u^*\ \
&\displaystyle\mbox{strongly in $L^2(\Omega)$ as
$\lambda\nearrow \lambda^*$},\\
&\displaystyle u_\lambda\to u^* \ \ &\displaystyle\mbox{a.e. in
$\Omega$ as $\lambda\nearrow \lambda^*$} .
\end{array} \right.
\end{equation}

Now we can proceed to obtain a contradiction. Multiplying by
$\varphi_1$ in $\eqref{Pla}$ and integrating over $\Omega$ we have
\begin{equation}\label{csapte} \displaystyle -\int_{\Omega} \varphi_1\, \Delta
u_{\lambda}=\lambda\int_{\Omega}f(u_{\lambda})
\varphi_1+\int_{\Omega} a(x)g(u_{\lambda})\varphi_1,\quad\mbox{
for all }\;0<\lambda<\lambda^*.
\end{equation}
On the other hand, by $(f1)$ it follows that $f(u_{\lambda})\geq
mu_{\lambda}$ in $\Omega,$ for all $0<\lambda<\lambda^*$.
Combining this with \eqref{csapte} we obtain
\begin{equation}\label{copt} \displaystyle
\lambda_1\int_{\Omega}u_{\lambda}\varphi_1\geq \lambda
m\int_{\Omega}u_{\lambda}\varphi_1+\int_{\Omega}
a(x)g(u_{\lambda})\varphi_1,\quad\mbox{ for all
}\;0<\lambda<\lambda^*.
\end{equation}
Notice that by $(g1)$, \eqref{car4} and the monotonicity of
$u_\lambda$ with respect to $\lambda$ we can apply the Lebesgue
convergence theorem to find
$$ \int_\Omega a(x)g(u_\lambda)\varphi_1\,dx\to \int_\Omega
a(x)g(u^*)\varphi_1\,dx\ \ \mbox{as } \lambda\nearrow \lambda_1.
$$ Passing to the limit in \eqref{copt} as $\lambda\nearrow\lambda^*,$ and
using \eqref{car4}, we obtain
$$
\displaystyle \lambda_1\int_{\Omega}u^*\varphi_1\geq
\lambda_1\int_{\Omega}u^*\varphi_1+\int_{\Omega}
a(x)g(u^*)\varphi_1. $$ Hence $\displaystyle
\int_{\Omega}a(x)g(u^*)\varphi_1=0,$ which is a contradiction.
Therefore $\displaystyle \lim_{\lambda\nearrow
\lambda^*}u_{\lambda}=+\infty$, uniformly on compact subsets of
$\Omega$. This concludes the proof.
\end{proof}

\section{Bifurcation and Asymptotics for a Singular Elliptic
Equation with Convection Term} Problems of this type arise in the
study of non-Newtonian fluids, boundary layer phenomena for
viscous fluids, chemical heterogeneous catalysts, cellular
automata and interacting particle systems with self--organized
criticality, as well as in the theory of Van der Waals
interactions in thin films spreading on solid surfaces (see, e.g.,
\cite{chayes,gennes,ock}).

 We are concerned in this section with singular elliptic
problems of the following type \begin{equation}\label{P}
\left\{\begin{tabular}{ll}
$-\Delta u=g(u)+\lambda|\nabla u|^p+\mu f(x,u)$ \quad & $\mbox{\rm in}\ \Omega,$\\
$u>0$ \quad & $\mbox{\rm in}\ \Omega,$\\
$u=0$ \quad & $\mbox{\rm on}\ \partial\Omega,$\\
\end{tabular} \right.
\end{equation}
where $\Omega\subset \RR^N$ $(N\geq 2)$ is a bounded domain with
smooth boundary, $0<p\leq 2$, and $\lambda$, $\mu\geq 0.$ We
suppose that
$f:\overline{\Omega}\times[0,\infty)\rightarrow[0,\infty)$ is a
H\"{o}lder continuous function which is non--decreasing with
respect to the second variable and is positive on
$\overline\Omega\times (0,\infty).$ We assume that
$g:(0,\infty)\rightarrow(0,\infty)$ is a H\"{o}lder continuous
function which is non--increasing and $\lim_{s\searrow
0}g(s)=+\infty.$ As in the previous section, we denote by
$\lambda_1$ the first eigenvalue of $(-\Delta)$ in
$H^1_0(\Omega).$ By the monotony of $g,$ there exists
$a:=\lim_{s\rightarrow\infty}g(s)\in[0,\infty).$

The next result concerns the case $\lambda=1$ and $1<p\leq 2.$

\begin{thm}\label{th11} Assume $\lambda=1$ and $1<p\leq 2.$ Then
the following properties hold true.

{\rm (i) } If $p=2$ and $a\geq \lambda_1,$ then problem \eqref{P}
has no solutions.

 {\rm (ii) } If either [$p=2$ and $a<\lambda_1$] or if $1<p<2,$ then
there exists $\mu^*>0$ such that problem \eqref{P} has at least
one classical solution for $\mu<\mu^*$ and no solutions exist if
$\mu>\mu^*.$
\end{thm}

In what follows the asymptotic behaviour of the nonlinear smooth
term $f(x,u)$ will play a decisive role. We impose the following
assumptions:

\noindent $\displaystyle (f4)\qquad$  there exists $c>0$ such that
$f(x,s)\geq cs$ for all $(x,s)\in \overline\Omega\times
[0,\infty);$

\noindent $\displaystyle (f5)\qquad$ the mapping $ (0,\infty)\ni
s\longmapsto f(x,s)/s$ is non--decreasing for all
$x\in\overline\Omega;$

\noindent $\displaystyle (f6)\qquad$ the mapping $ (0,\infty)\ni
s\longmapsto f(x,s)/s\quad\mbox{is non--increasing for all}\;\,
x\in\overline{\Omega};$

\noindent $\displaystyle (f7)\qquad$
$\lim_{s\rightarrow\infty}f(x,s)/s=0,\;\; \mbox{uniformly
for}\;\,x\in\overline{\Omega}.$

We first consider the case $\lambda=1$ and $0<p\leq 1$.

\begin{thm}\label{th21} Assume $\lambda=1$ and $0<p\leq 1.$ Then
the following properties hold true.

{\rm (i) } If $f$ satisfies either $(f4)$ or $(f5),$ then there
exists $\mu^*>0$ such that problem \eqref{P} has at least one
classical solution for $\mu<\mu^*$ and no solutions exist if
$\mu>\mu^*.$

 {\rm (ii) } If $0<p<1$ and $f$ satisfies
$(f6)-(f7),$ then problem \eqref{P} has at least one solution for
all $\mu\geq 0.$
\end{thm}

We now analyze the case $\mu=1$. Our framework is related to  the
sublinear case, described by assumptions $(f6)$ and $(f7)$.

\begin{thm}\label{th31} Assume $\mu=1$ and $f$ satisfies
assumptions $(f6)$ and $(f7).$ Then the following properties hold
true.

{\rm (i) } If $0<p<1,$ then problem \eqref{P} has at least one
classical solution for all $\lambda\geq 0$.

 {\rm (ii) } If $1\leq p\leq 2,$
then there exists $\lambda^*\in(0,\infty]$ such that problem
\eqref{P} has at least one classical solution for
$\lambda<\lambda^*$ and no solution exists if $\lambda>\lambda^*.$
Moreover, if $1<p\leq 2,$ then $\lambda^*$ is finite.
\end{thm}

Related to the above result we raise the following {\bf open
problem:} if $p=1$ and $\mu=1,$ is $\lambda^*$ a finite number?

Theorem \ref{th31} shows the importance of the convection term
$\lambda |\nabla u|^p$ in the singular problem \eqref{P}. Indeed,
according to Theorem \ref{th3} and for any $\mu>0$, the boundary
value problem
\begin{equation}\label{PQ} \left\{\begin{tabular}{ll}
$-\Delta u=u^{-\alpha}+\lambda|\nabla u|^p+\mu u^\beta$
\quad & $\mbox{\rm in}\ \Omega,$\\
$u>0$ \quad & $\mbox{\rm in}\ \Omega,$\\
$u=0$ \quad & $\mbox{\rm on}\ \partial\Omega$\\
\end{tabular} \right.
\end{equation}
has a unique solution, provided that $\lambda=0$ and $\alpha$,
$\beta\in (0,1)$. Theorem \ref{th31} shows that if $\lambda$ is
not necessarily 0, then the following situations may occur : (i)
problem \eqref{PQ} has solutions if $p\in (0,1)$ and for all
$\lambda\geq 0$; (ii) if $p\in (1,2)$ then there exists
$\lambda^*>0$ such that problem \eqref{PQ} has a solution for any
$\lambda<\lambda^*$ and no solution exists if $\lambda>\lambda^*.$

We give in what follows a complete description in the special case
$f\equiv 1$ and $p=2.$ More precisely, we consider the problem
\begin{equation}\label{Ppart}
\left\{\begin{tabular}{ll}
$-\Delta u=g(u)+\lambda|\nabla u|^2+\mu $ \quad & $\mbox{\rm in}\ \Omega,$\\
$u>0$ \quad & $\mbox{\rm in}\ \Omega,$\\
$u=0$ \quad & $\mbox{\rm on}\ \partial\Omega.$\\
\end{tabular} \right.
\end{equation}
A key role in this case will be played by the asymptotic behaviour
of the singular term $g$. In the statement of the next result we
remark some similarities with Theorem \ref{th5}.

\begin{thm}\label{th41}
The following properties hold true.

{\rm (i) } Problem \eqref{Ppart} has solution if and only if
$\lambda (a+\mu)<\lambda_1.$

 {\rm (ii) } Assume $\mu>0$ is fixed, $g$
is decreasing and let $ \lambda^*:=\lambda_1/(a+\mu).$ Then
problem \eqref{Ppart} has a unique solution $u_\lambda$ for all
$\lambda<\lambda^*$ and the sequence
$(u_\lambda)_{\lambda<\lambda^*}$ is increasing with respect to
$\lambda.$ Moreover, if $\limsup_{s\searrow 0}s^\alpha
g(s)<+\infty,$ for some $\alpha\in(0,1),$ then the sequence of
solutions $(u_\lambda)_{0<\lambda<\lambda^*}$ has the following
properties:

\qquad{\rm (ii1)} for all $0<\lambda<\lambda^*$ there exist two
positive constants $c_1$, $c_2$ depending on $\lambda$ such that
$c_1\,d(x)\leq u_\lambda\leq c_2\,d(x)$ in $\Omega;$

$\qquad${\rm (ii2)} $u_\lambda\in
C^{1,1-\alpha}(\overline\Omega)\cap C^2(\Omega);$

$\qquad${\rm (ii3)} $u_{\lambda}\longrightarrow +\infty$ as
$\lambda\nearrow \lambda^*$, uniformly on compact subsets of
$\Omega.$
\end{thm}

We refer to \cite{gr4} for complete proofs and further details.

\section{An Elliptic Problem with Strong Singular Nonlinearity and Convection Term}
We study the boundary value problem
\begin{equation}\label{prob1}
\left\{\begin{tabular}{ll}
$-\Delta u=p(x)g(u)+q(x)|\nabla u|^a$ \quad & $\mbox{\rm in}\ \Omega,$\\
$u>0$ \quad & $\mbox{\rm in}\ \Omega,$\\
$u=0$ \quad & $\mbox{\rm on}\ \partial\Omega,$\\
\end{tabular} \right.
\end{equation}
where $\Omega\subset\RR^N$ ($N\geq 2$) is a smooth bounded domain,
 $0<a<1$ and $q\in
C^{0,\alpha}(\overline\Omega)$, $q>0$ in $\overline\Omega$. The
potential $p\in C^{1}(\Omega)$ satisfies
\begin{equation}\label{condin}  c_1\,d(x)^\beta\leq |p(x)|\leq
c_2\,d(x)^\beta,\quad\mbox{ for all }\;x\in\Omega, \end{equation}
where $c_1$, $c_2>0,$  and $\beta$ is a real number. This
assumption shows that the potential $p(x)$ can admit a singular
boundary behaviour (corresponding to $\beta<0$).

Throughout this section we suppose  that $g\in C^{1}(0,\infty)$ is
a positive decreasing function such that $ \lim_{s\searrow
0}g(s)=+\infty$. The blow--up rate of $g$ at the origin is
described by the following assumption:

 \noindent $\displaystyle (g4)\qquad$ there exists $\gamma>\max\{1,\beta+1\}$ such that
$\lim_{s\searrow 0}s^\gamma g(s)\in(0,\infty).$

Observe that the stronger decay of the singular nonlinearity $g$
around the origin [described by our assumption $(g4)$] implies
that $g$ does {\bf not} obey the Keller--Osserman type condition
$(g3)$.

From \eqref{condin} we deduce that $p$ does not vanish in
$\Omega.$ Our first result concerns the case $p<0$ in $\Omega.$

\begin{thm}\label{th91} Assume that $g$ satisfies  $(g4),$ $p$ is negative
in $\Omega$, and condition \eqref{condin} is fulfilled. Then
problem \eqref{prob1} has no classical solutions.
\end{thm}

\begin{proof}
Let $\varphi_1$ be the normalized positive eigenfunction
corresponding to the first eigenvalue $\lambda_1$ of $(-\Delta)$
in $H^1_0(\Omega)$. Then $\lambda_1>0,$ $\varphi_1\in
C^2({\overline{\Omega}}),$ and
\begin{equation}\label{ineqphi}
\displaystyle C_1\,d(x)\leq \varphi_1(x)\leq C_2\,d(x),\quad
\mbox{  }\;x\in\Omega,
\end{equation}
for some positive constants $C_1$ and $C_2.$ From \eqref{condin}
and \eqref{ineqphi} it follows that there exist $\tau_1,\tau_2>0$
such that
\begin{equation}\label{ipotezak}
\tau_1\varphi_1(x)^\beta\leq |p(x)|\leq
\tau_2\varphi_1(x)^\beta,\quad\mbox{ for all }\;x\in\Omega.
\end{equation}

Fix $C>0$ such that $\|q\|_{\infty}^2C^{a-1}<\lambda_1$ and define
$\psi:[0,\infty)\rightarrow[0,\infty)$ by $\psi(s)=s^a/(s^2+C).$
Then $\psi$ attains its maximum at $\bar
s=\left[Ca/(2-a)\right]^{1/2}.$ Hence
$$\displaystyle \psi(s)\leq \psi(\bar s)=\frac{a^{a/2}(2-a)^{(2-a)/2}}{2C^{1-a/2}},
\quad\mbox{ for all }\;s\geq 0.$$ An elementary computation shows
that
\begin{equation}\label{inq} \displaystyle s^a\leq
C^{a/2-1}s^2+C^{a/2},\qquad \mbox{ for all }\;s\geq 0.
\end{equation}

Arguing by contradiction, let us assume that problem \eqref{prob1}
has a classical solution $U.$ Consider the perturbed problem
\begin{equation}\label{pert}
\left\{\begin{tabular}{ll}
$\displaystyle -\Delta u=p(x)g(u+\varepsilon)+A|\nabla u|^2+B$ \quad & $\mbox{\rm in}\ \Omega,$\\
$u>0$ \quad & $\mbox{\rm in}\ \Omega,$\\
$u=0$ \quad & $\mbox{\rm on}\ \partial\Omega,$\\
\end{tabular} \right.
\end{equation}
where $\varepsilon>0,$ and $A=\|q\|_{\infty} C^{a/2-1},$
$B=\|q\|_{\infty} C^{a/2}.$ By virtue of \eqref{inq} it follows
that $U$ is a sub-solution of \eqref{pert}. Set $v=e^{Au}-1.$ Then
problem \eqref{pert} becomes
\begin{equation}\label{pertv}
\left\{\begin{tabular}{ll} $\displaystyle -\Delta
v=Ap(x)(v+1)g\Big(\frac{1}{A}\ln(v+1)+\varepsilon\Big)+AB(v+1)$
\quad &
$\mbox{\rm in}\ \Omega,$\\
$v>0$ \quad & $\mbox{\rm in}\ \Omega,$\\
$v=0$ \quad & $\mbox{\rm on}\ \partial\Omega.$\\
\end{tabular} \right.
\end{equation}
 We first remark that $V=e^{AU}-1$ is a sub-solution
of \eqref{pertv}. On the other hand, since $AB<\lambda_1,$ we
conclude that there exists $w\in C^2(\overline\Omega)$ such that
\begin{equation}\label{variab}
\left\{\begin{tabular}{ll}
$\displaystyle -\Delta w=AB(w+1)$ \quad & $\mbox{\rm in}\ \Omega,$\\
$w>0$ \quad & $\mbox{\rm in}\ \Omega,$\\
$w=0$ \quad & $\mbox{\rm on}\ \partial\Omega.$\\
\end{tabular} \right.
\end{equation}
Moreover, the maximum principle yields
\begin{equation}\label{estim}
c_1\varphi_1\leq w\leq c_2\varphi_2\quad\mbox{ in }\;\Omega,
\end{equation}
for some positive constants $c_1$ and $c_2>0.$ It is clear that
$w$ is a super-solution of \eqref{pertv}. We claim that $V\leq w.$
To this aim, it suffices to prove that $U\leq W$ in $\Omega,$
where $W=A^{-1}\ln(w+1)$ verifies
$$\left\{\begin{tabular}{ll}
$\displaystyle -\Delta W=A|\nabla W|^2+B$ \quad & $\mbox{\rm in}\ \Omega,$\\
$W>0$ \quad & $\mbox{\rm in}\ \Omega,$\\
$W=0$ \quad & $\mbox{\rm on}\ \partial\Omega.$\\
\end{tabular} \right.$$

Assuming the contrary, we get that $\max_{x\in\overline\Omega}
(U-W)>0$ is achieved in some point $x_0\in\Omega.$ Then
$\nabla(U-W)(x_0)=0$ and
$$
\displaystyle 0\leq -\Delta(U-W)(x_0)
=p(x_0)g(U(x_0))+q(x_0)|\nabla U|^a(x_0)-A|\nabla W|^2(x_0)-B<0,
$$ which is a contradiction. Hence $U\leq W$ in $\Omega,$ that is,
$V\leq w$ in $\Omega.$ By the sub and super-solution method we
deduce that there exists $v_\varepsilon\in C^2(\overline\Omega)$ a
solution of problem \eqref{pertv} such that
\begin{equation}\label{pertmarg}
V\leq v_\varepsilon\leq w\quad\mbox{ in }\;\Omega.
\end{equation}

Now we proceed to get our contradiction. Integrating in
\eqref{pertv} and taking into account the fact that $p$ is
negative, we deduce
$$\displaystyle -\int_\Omega \Delta v_\varepsilon dx-A\int_\Omega
p(x)g\Big(\frac{1}{A}\ln(v_\varepsilon+1)+\varepsilon\Big)dx\leq
AB\int_\Omega (w+1)dx.$$ Using monotonicity of $g$ and the fact
that $\ln(v_\varepsilon+1)\leq v_\varepsilon$ in $\Omega,$ the
above inequality yields
$$\displaystyle -\int_{\partial\Omega} \frac{\partial v_\varepsilon}{\partial n}ds-A\int_\Omega
p(x)g\Big(\frac{v_\varepsilon}{A}+\varepsilon\Big)dx\leq
AB(\|w\|_\infty+1)|\Omega|<+\infty.$$ Since $\partial
v_{\varepsilon}/\partial n\leq 0$ on $\partial\Omega,$ the above
relation implies
\begin{equation}\label{raport}
\displaystyle
-\int_{\Omega}p(x)g\Big(\frac{v_{\varepsilon}}{A}+\varepsilon\Big)dx\leq
M,
\end{equation}
where $M=B(\|w\|_\infty+1)|\Omega|.$ Now, relations
\eqref{pertmarg} and \eqref{raport} imply $$0\leq
-\int_{\Omega}p(x)g\left(\frac{w}{A}+\varepsilon\right)dx\leq M.$$
Therefore, for any compact subset $\omega\subset\subset\Omega$ we
have
$$\displaystyle
0\leq -\int_{\omega}p(x)g\Big(\frac{w}{A}+\varepsilon\Big)dx\leq
M.$$ Passing to the limit as $\varepsilon\searrow 0$ in the above
inequality, it follows that
$$-\int_{\omega}p(x)g\left(\frac{w}{A}\right)dx\leq M,\qquad\mbox{for all
 $\omega\subset\subset\Omega$}.$$ This yields
\begin{equation}\label{fin}
\displaystyle -\int_{\Omega}p(x)g\Big(\frac{w}{A}\Big)dx\leq
M.\end{equation} On the other hand, the hypothesis $(g4)$ combined
with \eqref{estim} implies $g\big(w/A\big)\geq
c_0\varphi_1^{-\gamma}$ in $\Omega,$ for some $c_0>0.$ The last
inequality together with \eqref{ipotezak} and \eqref{fin} produces
$c\int_{\Omega}\varphi_1^{\beta-\gamma}dx\leq M,$ where
$\beta-\gamma<-1.$ But, by a result of Lazer and McKenna (see
\cite{lm1}), $\int_{\Omega}\varphi_1^{-s}dx<+\infty$ if and only
if $s<1.$ This contradiction shows that problem \eqref{prob1} has
no classical solutions and the proof  is now complete.
\end{proof}

The situation changes radically in the case where $p$ is positive
in $\Omega$, as established in the next result.

\begin{thm}\label{th92} Assume that $g$ satisfies $(g4)$
 and the potential
$p(x)$ is positive and fulfills \eqref{condin}. Then the following
properties hold true.

 {\rm(i) } If $\beta\leq -2,$ then problem
\eqref{prob1} has no classical solutions.

 {\rm(ii) } If
$\beta>-2$, then problem \eqref{prob1} has a unique solution $u$
which, moreover, has the following properties:

{\rm(ii1) } there exist $M$, $m>0$ such that $$
m\,d(x)^{(2+\beta)/(1+\gamma)} \leq u(x)\leq
M\,d(x)^{(2+\beta)/(1+\gamma)},\quad\mbox{ for all }\;x\in\Omega;
$$

{\rm(ii2) } if $\beta\geq\max\{0,\gamma-3\},$ then $u$ is in
$H^1_0(\Omega);$

{\rm(ii3) } if $2\beta\leq\gamma-3,$ then $u$ does not belong to
$H^1_0(\Omega).$
\end{thm}

We refer to \cite{gr5} for the proof of Theorem \ref{th92}, as
well as for a result concerning the entire solutions of problem
\eqref{prob1}.

\subsection*{Acknowledgment}  I am greatly indebted to
Professor Haim Brezis, for his highest level guidance during my
PhD and Habilitation theses at the Universit\'e Pierre et Marie
Curie (Paris~6), as well as for suggesting to me several modern
research subjects and directions of interest. {\it Bonne
anniversaire, mon Professeur!}


\begin{thebibliography}{1}
\bibitem{bbc} P. B\'enilan, H. Brezis, and M. Crandall,
\textit{A semilinear equation in $L^1(\RR^N)$.}  Ann. Scuola Norm.
Sup. Pisa \textbf{4} (1975), 523--555.

\bibitem{chayes} J. T. Chayes, S. J. Osher, and J. V. Ralston, \textit{On singular diffusion
equations with applications to self-organized criticality.} Comm.
Pure Appl. Math. \textbf{46} (1993), 1363--1377.

\bibitem{cgr} F.-C. C\^{\i}rstea, M. Ghergu, and V. R\u adulescu, \textit{Combined effects
of asymptotically linear and singular nonlinearities in
bifurcation problems of Lane-Emden-Fowler type.} J.~Math. Pures
Appl., in press.

\bibitem{cp} M. Coclite and G. Palmieri, \textit{On a singular nonlinear Dirichlet
problem.} Commun. Partial Diff. Equations \textbf{14} (1989),
1315--1327.

\bibitem{crt} M. G. Crandall, P. H. Rabinowitz,
and L. Tartar, \textit{On a Dirichlet problem with a singular
nonlinearity.} Commun. Partial Diff. Equations \textbf{2} (1977),
193--222.

\bibitem{gennes} P. G. de Gennes, \textit{Wetting: statics and dynamics.} Review of
Modern Physics \textbf{57} (1985), 827x-863.

\bibitem{gr1} M. Ghergu and V. R\u adulescu, \textit{Bifurcation and asymptotics for the
Lane-Emden-Fowler equation.} C. R. Acad. Sci. Paris, Ser.~I
 \textbf{337} (2003), 259--264.

\bibitem{gr2} M. Ghergu and V. R\u adulescu, \textit{Sublinear singular
elliptic problems with two parameters.} J.
Differential Equations
 \textbf{195} (2003), 520--536.

\bibitem{gr3} M. Ghergu and V. R\u adulescu, \textit{Bifurcation for a class of singular elliptic
problems with quadratic convection term.} C. R. Acad. Sci. Paris,
Ser.~I \textbf{338} (2004), 831--836.

\bibitem{gr4} M. Ghergu and V. R\u adulescu, \textit{Multiparameter
bifurcation and asymptotics for the singular Lane-Emden-Fowler
equation with a convection term.} Proc. Royal Soc. Edinburgh Sect.
A, in press.

\bibitem{gr5} M. Ghergu and V. R\u adulescu, \textit{Singular elliptic problems with
sublinear convection term and Kato potential in anisotropic
media}, in preparation.

\bibitem{h} L. H\"ormander, \textit{The Analysis of Linear
Partial Differential Operators I.} Springer, Berlin, 1983.

\bibitem{lm1} A. C. Lazer and P. J. McKenna, \textit{On a
singular nonlinear elliptic boundary value problem.}  Proc. Amer.
Math. Soc. \textbf{3} (1991), 720--730.

\bibitem{mr} P. Mironescu and V. R\u adulescu, \textit{The study
of a bifurcation problem associated to an asymptotically linear
function.}  Nonlinear Anal., T.M.A. \textbf{26} (1996), 857--875.

\bibitem{ock} J. Ockendon, S. Howison, A. Lacey, and A. Movchan,
\textit{Applied Partial Differential Equations}, Oxford University
Press, 2003.

\bibitem{shi} J. Shi and M. Yao, \textit{On a singular
nonlinear semilinear elliptic problem.} Proc. Roy. Soc. Edinburgh
Sect. A \textbf{128} (1998), 1389--1401.


\end{thebibliography}
\end{document}